\magnification=1200
\overfullrule=0pt
{\centerline {\bf A general variational principle and some of its 
applications}}\par
\bigskip
\bigskip
{\centerline {BIAGIO RICCERI}}\par
\bigskip
\bigskip
\bigskip
{\bf 1. Introduction}\par
\bigskip
The aim of this paper is to establish Theorem 2.1 below and, mainly,
to show how one can derive from it, in an absolutely transparent way,
a series of consequences about the local minima of a functional
of the type $\Phi+\lambda\Psi$, where $\lambda$ is a positive real
number, and, in concrete situations, $\Phi$ and $\Psi$ are sequentially
weakly lower semicontinuous functionals defined on a subset of a reflexive
Banach space.\par
\smallskip
 Ultimately, our end is to apply, via the variational
methods, the basic theory developed below
to differential equations. In this connection,
Theorem 2.5 can be regarded as the main result of this paper. Specific
applications of Theorem 2.5 to differential equations will systematically
be presented in a series of successive papers. Here we limit ourselves
to give a sample of application to a class of elliptic equations involving 
the critical Sobolev exponent (Theorem 2.8).\par
\smallskip
 We also derive from Theorem 2.5 a result
on fixed points of potential operators in Hilbert spaces (Theorem 2.7)
of which the following is a corollary:\par
\medskip
THEOREM A. - {\it Let $X$ be a real Hilbert space, and let
$A:X\to X$ be a potential operator, with
 sequentially weakly upper semicontinuous 
potential $P$ satisfying 
$$\liminf_{r\to +\infty}{{\sup_{\|x\|\leq r}P(x)}\over {r^2}}<{{1}\over
{2}}<
\limsup_{r\to +\infty}{{\sup_{\|x\|\leq r}P(x)}\over {r^2}}\ .$$
Then, the set of all fixed points of $A$
is unbounded.}\par
\bigskip
{\bf 2. Results}\par
\bigskip
Our abstract basic result is as follows:\par
\medskip
THEOREM 2.1. - {\it Let $X$ be a topological space, and let
$\Phi, \Psi : X\to {\bf R}$ be two sequentially lower semicontinuous
functions. Denote by $I$ the set of all $\rho>\inf_{X}\Psi$ such that
the set  $\Psi^{-1}(]-\infty,\rho[)$ is contained in
some sequentially compact subset of $X$. Assume that $I\neq\emptyset$.
For each $\rho\in I$, denote by ${\cal F}_{\rho}$ the family of
all sequentially compact subsets of $X$ containing 
$\Psi^{-1}(]-\infty,\rho[)$, and put
$$\alpha(\rho)=\sup_{K\in {\cal F}_{\rho}}\inf_{K}\Phi.$$
Then, for each $\rho\in I$ and each $\lambda$ satisfying
$$\lambda>\inf_{x\in \Psi^{-1}(]-\infty,\rho[)}{{\Phi(x)-\alpha(\rho)}
\over {\rho-\Psi(x)}}$$
the restriction of the function $\Phi+\lambda\Psi$ to 
$\Psi^{-1}(]-\infty,\rho[)$ has a global minimum.}\par
\smallskip
PROOF. Fix $\rho\in I$ and $\lambda$ as in the conclusion. Observe that
$$\inf_{x\in \Psi^{-1}(]-\infty,\rho[)}{{\Phi(x)-\alpha(\rho)}
\over {\rho-\Psi(x)}}=\inf_{r>-\alpha(\rho)}
\inf_{x\in \Psi^{-1}(]-\infty,\rho[)}{{\Phi(x)+r}
\over {\rho-\Psi(x)}}\ .$$
Consequently, we can fix $r^{*}>-\alpha(\rho)$ so that
$$\inf_{x\in \Psi^{-1}(]-\infty,\rho[)}{{\Phi(x)+r^{*}}
\over {\rho-\Psi(x)}}<\lambda. \eqno{(1)}$$
On the other hand, since $-r^{*}<\alpha(\rho)$, there is a sequentially
compact subset $K$ of $X$ containing $\Psi^{-1}(]-\infty,\rho[)$
such that
$$-r^{*}<\inf_{K}\Phi.$$
For each $r>-\inf_{K}\Phi$, put
$$\beta_{\rho}(r)=\sup_{x\in \Psi^{-1}(]-\infty,\rho[)}
{{\Phi(x)+r}\over {\Psi(x)-\rho}}\ .$$
Observe that the function ${{\Phi(\cdot)+r}\over {\Psi(\cdot)-\rho}}$
is negative in $\Psi^{-1}(]-\infty,\rho[)$. We now show that it there
attains the value $\beta_{\rho}(r)$. To this end, fix a sequence
$\{x_{n}\}$ in $\Psi^{-1}(]-\infty,\rho[)$ such that
$$\lim_{n\to \infty}{{\Phi(x_{n})+r}\over {\Psi(x_{n})-\rho}}
=\beta_{\rho}(r). \eqno{(2)}$$
Since $K$ is sequentially compact, $\{x_{n}\}$
admits a subsequence $\{x_{n_{k}}\}$ converging to a point $x^{*}\in K$.
Put
$$l=\liminf_{k\to \infty}\Psi(x_{n_{k}}).$$
We claim that $l<\rho$. Indeed, if it was $l=\rho$, since
$$l\leq \limsup_{k\to \infty}\Psi(x_{n_{k}})\leq \rho,$$
we would have
$$\lim_{k\to \infty}\Psi(x_{n_{k}})=\rho,$$
and so, from $(2)$, we would infer
$$\lim_{k\to \infty}(\Phi(x_{n_{k}})+r)=0.$$
Then, since $\Phi$ is sequentially lower semicontinuous, it would follow
$$\Phi(x^{*})+r\leq 0,$$
from which
$$r\leq -\inf_{K}\Phi\ ,$$
against the choice of $r$. Since $\Psi$ is sequentially lower 
semicontinuous, we have
$$\Psi(x^*)\leq l<\rho. \eqno{(3)}$$
Now, extract from $\{x_{n_{k}}\}$ a subsequence $\{x_{n_{k_{p}}}\}$
so that
$$l=\lim_{p\to \infty}\Psi(x_{n_{k_{p}}}).$$
>From $(2)$, we then get
$$\lim_{p\to \infty}(\Phi(x_{n_{k_{p}}})+r)=\beta_{\rho}(r)(l-\rho).$$
So, by the sequential lower semicontinuity of $\Phi$, we have
$$0<\Phi(x^{*})+r\leq \beta_{\rho}(r)(l-\rho). \eqno{(4)}$$
Now, from $(3)$ and $(4)$, we directly obtain
$$\beta_{\rho}(r)\leq {{\Phi(x^{*})+r}\over {\Psi(x^{*})-\rho}}$$
which shows our claim. Clearly, the function $\beta_{\rho}$, as
the supremum of affine functions, is convex in the open interval
$]-\inf_{K}\Phi,+\infty[$. Hence, it is continuous there.
Also, observe that it is not bounded below. Indeed,
one has
$$\beta_{\rho}(r)\leq {{r+\inf_{K}\Phi}\over
{\inf_{K}\Psi-\rho}}$$
for all $r>-\inf_{K}\Phi$. 
Therefore, recalling $(1)$, since $r^{*}>-\inf_{K}\Phi$, there
exists $r_{0}>-\inf_{K}\Phi$ such that $$\beta_{\rho}(r_{0})=
-\lambda\ .$$
Finally, let $x_{0}\in \Psi^{-1}(]-\infty,\rho[)$ be such that
$$\beta_{\rho}(r_{0})={{\Phi(x_{0})+r_{0}}\over {\Psi(x_{0})-\rho}}\ .$$
Hence, for each $x\in \Psi^{-1}(]-\infty,\rho[)$, one has
$${{\Phi(x)+r_{0}}\over {\Psi(x)-\rho}}\leq -\lambda$$
and so
$$\Phi(x)+r_{0}+\lambda(\Psi(x)-\rho)\geq 0=\Phi(x_{0})+r_{0}+\lambda
(\Psi(x_{0})-\rho)$$
that gives
$$\Phi(x)+\lambda\Psi(x)\geq \Phi(x_{0})+\lambda\Psi(x_{0})\ .$$
This concludes the proof.\hfill $\bigtriangleup$
\medskip
REMARK 2.1. - Concerning the statement of Theorem 2.1, observe that 
when the space $X$ is Hausdorff, the set $H=\cap\{K : K\in
{\cal F}_{\rho}\}$ belongs to ${\cal F}_{\rho}$, and so we
have
$$\alpha(\rho)=\inf_{H}\Phi.$$
\smallskip
One of the most significant features of Theorem 2.1 is the possibility
to get from it multiplicity results for local minima. Indeed,
we have\par
\medskip
THEOREM 2.2. - {\it Let the assumptions of Theorem 2.1 be satisfied.
In addition, suppose
$$\sup I=+\infty$$
and 
$$\gamma<+\infty\ ,$$
where
$$\gamma=\liminf_{\rho\to +\infty}
\inf_{x\in \Psi^{-1}(]-\infty,\rho[)}{{\Phi(x)-\alpha(\rho)}
\over {\rho-\Psi(x)}}\ .$$
Finally, denote by $\tau$ the weakest topology on $X$ for which
$\Psi$ is upper semicontinuous.\par
Then, for each $\lambda>\gamma$, the following alternative
holds: either $\Phi+\lambda \Psi$ has a global
minimum, or there exists a sequence $\{x_{n}\}$ of $\tau$-local minima
of $\Phi+\lambda \Psi$ such that
$$\lim_{n\to \infty}\Psi(x_{n})=+\infty\ .$$}
PROOF. Let $\lambda>\gamma$. Then, we can fix a sequence $\{\rho_{n}\}$ in
$I$, with $\lim_{n\to \infty}\rho_{n}=+\infty$, such that
$$\inf_{x\in \Psi^{-1}(]-\infty,\rho_{n}[)}{{\Phi(x)-\alpha(\rho_{n})}
\over {\rho_{n}-\Psi(x)}}<\lambda$$
for all $n\in {\bf N}$. Consequently, thanks to Theorem 2.1, for each
$n\in {\bf N}$, there is $x_{n}\in \Psi^{-1}(]-\infty,\rho_{n}[)$
such that 
$$\Phi(x_{n})+\lambda \Psi(x_{n})\leq \Phi(x)+\lambda \Psi(x) \eqno{(5)}$$
for all $x\in \Psi^{-1}(]-\infty,\rho_{n}[)$. Now, if 
$\lim_{n\to \infty}\Psi(x_{n})=+\infty$, we are done, since
each set $\Psi^{-1}(]-\infty,\rho_{n}[)$ is $\tau$-open. Thus,
suppose that $\liminf_{n\to \infty}\Psi(x_{n})<+\infty$. Hence,
 there are an increasing sequence $\{n_{k}\}$ in
${\bf N}$ and a constant $c\in I$ such that
$$\Psi(x_{n_{k}})<c$$
for all $k\in {\bf N}$. But the set $\Psi^{-1}(]-\infty,c[)$ is contained
in some sequentially compact subset of $X$, and so there is a
subsequence $\{x_{n_{k_{r}}}\}$ converging to some $x^{*}\in X$. 
Finally, fix $x\in X$. Since we definitively have $\rho_{n_{k_{r}}}
>\Psi(x)$, taking into account the sequential lower semicontinuity
of $\Phi+\lambda \Psi$, from $(5)$ we get
$$\Phi(x^{*})+\lambda \Psi(x^{*})\leq \liminf_{r\to +\infty}
(\Phi(x_{n_{k_{r}}})+\lambda \Psi(x_{n_{k_{r}}}))\leq
 \Phi(x)+\lambda \Psi(x)\ .$$
Hence, $x^*$ is a global minimum of $\Phi+\lambda \Psi$.
\hfill $\bigtriangleup$
\medskip
THEOREM 2.3. - {\it Let the assumptions of Theorem 2.1 be satisfied.
In addition, suppose
$$\delta<+\infty\ ,$$
where
$$\delta=\liminf_{\rho\to \left ( {\inf\limits _{X}\Psi}\right ) ^{+}}
\inf_{x\in \Psi^{-1}(]-\infty,\rho[)}{{\Phi(x)-\alpha(\rho)}
\over {\rho-\Psi(x)}}\ .$$
Finally, denote by $\tau$ the weakest topology on $X$ for which
$\Psi$ is upper semicontinuous.\par
Then, for each $\lambda>\delta$, there exists a sequence 
of $\tau$-local minima of $\Phi+\lambda \Psi$ which
converges to a global minimum of $\Psi$.}\par
\smallskip
PROOF. Let $\lambda>\delta$. Fix a sequence $\{\rho_{n}\}$ in $I$,
with $\lim_{n\to \infty}\rho_{n}=\inf_{X}\Psi$, such that
$$\inf_{x\in \Psi^{-1}(]-\infty,\rho_{n}[)}{{\Phi(x)-\alpha(\rho_{n})}
\over {\rho_{n}-\Psi(x)}}<\lambda$$
for all $n\in {\bf N}$. Thanks to Theorem 2.1, for each $n\in {\bf N}$,
the restriction of $\Phi+\lambda \Psi$ to $\Psi^{-1}(]-\infty,\rho_{n}[)$
has a global minimum, say $x_{n}$. Hence, $x_{n}$ is a $\tau$-local
minimum of $\Phi+\lambda \Psi$. Finally, the sequence $\{x_{n}\}$
lies in $\Psi^{-1}(]-\infty,\max_{n\in {\bf N}}\rho_{n}[)$ (which is,
in turn, contained in a sequentially compact subset of $X$), and so
it admits a subsequence converging to a point which, by the
sequential lower semicontinuity of $\Psi$, is a global minimum
of $\Psi$.\hfill $\bigtriangleup$\par
\medskip
We now derive from the previous abstract theorems more concrete
results in reflexive Banach spaces.
\smallskip
THEOREM 2.4. - {\it Let $E$ be a reflexive real Banach space, $X$ a closed,
convex, unbounded subset of $E$, and $\Phi, \Psi :X\to {\bf R}$ two convex
functionals, with $\Phi$ lower semicontinuous and $\Psi$ continuous and
satisfying $\lim_{x\in X, \|x\|\to +\infty}\Psi(x)=+\infty$. Put
$$\lambda^{*}=\inf_{\rho>\inf\limits _{X}\Psi}\inf_{x\in
\Psi^{-1}(]-\infty,\rho[)}
{{\Phi(x)-\inf_{\Psi^{-1}(]-\infty,\rho])}\Phi}\over
 {\rho-\Psi(x)}}\ .$$ 
Then, for each $\lambda>\lambda^{*}$, the functional $\Phi+\lambda
\Psi$ has a global minimum in $X$. Moreover, if $\lambda^{*}>0$,
for each $\mu<\lambda^{*}$, the functional $\Phi+\mu\Psi$
has no global minima in $X$.}\par
\smallskip
PROOF. It is clear that we can apply Theorem 2.1 endowing $X$ with
the relativization of the weak topology. In particular, for
each $\rho>\inf_{X}\Psi$, the set $\Psi^{-1}(]-\infty,
\rho])$ is sequentially weakly compact owing to the reflexivity of
$E$ and to the coercivity of $\Psi$. Moreover, by the convexity of
$\Psi$, the same set is equal to the closure of $\Psi^{-1}(]-\infty,\rho[)$.
>From this, it follows that $\Psi^{-1}(]-\infty,\rho])$ is the
smallest sequentially weakly compact subset of $X$ containing
$\Psi^{-1}(]-\infty,\rho[)$. Hence, by Remark 2.1, we have
$$\alpha(\rho)=\inf_{\Psi^{-1}(]-\infty,\rho])}\Phi\ .$$
Now, let $\lambda>\lambda^{*}$ and choose $\rho>\inf_{X}\Psi$ so that
$$\lambda>\inf_{x\in
\Psi^{-1}(]-\infty,\rho[)}
{{\Phi(x)-\inf_{\Psi^{-1}(]-\infty,\rho])}\Phi}\over
 {\rho-\Psi(x)}}\ .$$ 
Then, Theorem 2.1 ensures that the restriction of the functional
$\Phi+\lambda\Psi$ to $\Psi^{-1}(]-\infty,\rho[)$ has a global minimum,
say $x_{0}$.
But, by assumption, $\Psi$ is (strongly) continuous, and so
the set $\Psi^{-1}(]-\infty,\rho[)$ is (strongly) open in $X$. In
other words, $x_{0}$ is a local minimum for $\Phi+\lambda\Psi$ in
the strong topology. Since $\Phi+\lambda\Psi$ is convex, $x_{0}$
is actually a global minimum for $\Phi+\lambda\Psi$ in $X$.\par
Now, assume that $\lambda^{*}>0$.
 Let $\mu$ be such that the functional $\Phi+\mu\Psi$ has
a global minimum in $X$, say $x_{1}$. We claim that
 $\mu\geq \lambda^{*}$ from
which, of course, the second part of the conclusion follows. 
Indeed, fix $\rho>\Psi(x_{1})$ and choose $x_{2}\in
\Psi^{-1}(]-\infty,\rho])$ so that
$$\Phi(x_{2})=\inf_{\Psi^{-1}(]-\infty,\rho])}\Phi\ .$$
Clearly, since $\lambda^{*}>0$, we have $\Psi(x_{2})=\rho$.
Hence, we get
$$\mu\geq {{\Phi(x_{1})-\Phi(x_{2})}\over {\Psi(x_{2})-\Psi(x_{1})}}=
{{\Phi(x_{1})-\alpha(\rho)}\over {\rho-\Psi(x_{1})}}\geq
\inf_{x\in \Psi^{-1}(]-\infty,\rho[)}
{{\Phi(x)-\alpha(\rho)}\over {\rho-\Psi(x)}}\geq \lambda^{*}\ ,$$
as claimed.\hfill $\bigtriangleup$\par
\medskip
The next result groups together the 
versions of Theorems 2.1, 2.2 and 2.3 which are directly applicable
to differential equations.\par
\medskip
THEOREM 2.5. - {\it Let $X$ be a reflexive real Banach space,
and let $\Phi, \Psi
:X\to {\bf R}$ be two sequentially weakly lower semicontinuous and
G\^ateaux differentiable functionals. Assume also that $\Psi$ is
(strongly) continuous and satisfies $\lim_{\|x\|\to +\infty}\Psi(x)=
+\infty$. For each $\rho>\inf_{X}\Psi$, put 
$$\varphi(\rho)=\inf_{x\in \Psi^{-1}(]-\infty,\rho[)}
{{\Phi(x)-\inf_{\overline {(\Psi^{-1}(]-\infty,\rho[))}_{w}}\Phi}
\over {\rho-\Psi(x)}}\ ,$$
where $\overline {(\Psi^{-1}(]-\infty,\rho[))}_{w}$ is the closure
of $\Psi^{-1}(]-\infty,\rho[)$ in the weak topology. Furthermore,
set
$$\gamma=\liminf_{\rho\to +\infty}\varphi(\rho)$$
and
$$\delta=\liminf_{\rho\to 
\left ( {\inf\limits _{X}\Psi}\right ) ^{+}}\varphi(\rho)\ .$$
Then, the following conclusions hold:\par
\noindent
$(a)$\hskip 10pt For each $\rho>\inf_{X}\Psi$ and each $\lambda>
\varphi(\rho)$, the functional $\Phi+\lambda\Psi$ has a critical
point which lies in $\Psi^{-1}(]-\infty,\rho[)$.\par
\noindent
$(b)$\hskip 10pt If $\gamma<+\infty$, then, for each $\lambda>
\gamma$, the following alternative holds: either $\Phi+\lambda
\Psi$ has a global minimum, or the set of all critical points
of $\Phi+\lambda \Psi$ is unbounded.\par
\noindent
$(c)$\hskip 10pt If $\delta<+\infty$, then, for each $\lambda>
\delta$, the following alternative holds:
either there exists a global minimum of $\Psi$ which is
a local minimum of $\Phi+\lambda\Psi$, or there exists a
sequence of pairwise distinct critical points of
$\Phi+\lambda\Psi$ which weakly converges to a global minimum of
$\Psi$.}\par
\smallskip
PROOF. Endow $X$ with the weak topology.
As in the proof of Theorem 2.4, it is seen that 
$\overline {(\Psi^{-1}(]-\infty,\rho[))}_{w}$ is the smallest
sequentially weakly compact subset of $X$ containing
$\Psi^{-1}(]-\infty,\rho[)$. So, by Remark 2.1, we have
$$\alpha(\rho)=\inf_{\overline {(\Psi^{-1}(]-\infty,\rho[))}_{w}}\Phi
\ .$$
Now, $(a)$ follows at once from Theorem 2.1 since 
any global minimum of the restriction of $\Phi+\lambda \Psi$ to
$\Psi^{-1}(]-\infty,\rho[)$ is, by the continuity of $\Psi$, a local minimum
of $\Phi+\lambda \Psi$ in the strong topology.\par
Assume $\gamma<+\infty$, and let $\lambda>\gamma$. Then, by Theorem 2.2, 
 if $\Phi+\lambda\Psi$
has no global minima, there exists a sequence $\{x_{n}\}$ of
local minima of $\Phi+\lambda\Psi$ in the strong
topology such that $\lim_{n\to \infty}
\Psi(x_{n})=+\infty$. Owing to the coercivity of $\Psi$, the
sequence $\{x_{n}\}$ is clearly unbounded, and so $(b)$ follows.\par
Finally, assume $\delta<+\infty$, and 
let $\lambda>\delta$. Assume that no global minimum of
$\Psi$ is a local minimum of $\Phi+\lambda\Psi$. By Theorem
2.3, there exists a sequence of local minima of $\Phi+\lambda\Psi$ in
the strong topology
which weakly converges to a global minimum of $\Psi$. But then,
since $X$ equipped
with the weak topology is Hausdorff, there is a subsequence whose
terms are pairwise distinct, and we are done.\hfill
$\bigtriangleup$\par 
\medskip
It is also worth noticing the following corollary of
Theorem 2.5.\par
\medskip
THEOREM 2.6. - {\it Let the assumptions of Theorem 2.5 be satisfied.
Assume also that $\delta<+\infty$.
\par
Then, either the system
$$\cases {\Phi'(x)=0 \cr & \cr \Psi'(x)=0\cr}$$
has at least one solution, or, for each $\lambda>\delta$,
there exists a
sequence of pairwise distinct critical points of
$\Phi+\lambda\Psi$ which weakly converges to a global minimum of
$\Psi$. In particular, the first case of the alternative
occurs when there exists a Hausdorff
vector topology on $X^*$ with respect to 
which the operators $\Phi'$ and $\Psi'$ are sequentially weakly 
continuous.}\par
\smallskip
PROOF. Assume that the system 
$$\cases {\Phi'(x)=0 \cr & \cr \Psi'(x)=0\cr}$$
has no solutions. Let $\lambda>\delta$. Then, of course, no
global minimum of $\Psi$ can also be a local minimum of $\Phi+\lambda\Psi$.
Consequently, by Theorem 2.5, there is a sequence $\{x_{n}\}$ of
pairwise distinct critical points of $\Phi+\lambda\Psi$ which
weakly converges to a global minimum of $\Psi$, say $x^*$. So,
we have
$$\Phi'(x_{n})+\lambda\Psi'(x_{n})=0 \eqno{(6)}$$
for all $n\in {\bf N}$.
 Finally, observe that there is no Hausdorff vector topology on
$X^*$ with respect to which the
operators $\Phi'$ and $\Psi'$ are  
sequentially weakly continuous, since, otherwise, passing to the
limit in $(6)$, we would have
$$\Phi'(x^{*})+\lambda\Psi'(x^{*})=0$$
from which
$$\cases {\Phi'(x^{*})=0 \cr & \cr \Psi'(x^{*})=0\ ,\cr}$$
against our assumption.\hfill $\bigtriangleup$
\medskip
We now point out another corollary of Theorem 2.5 concerning fixed points of
potential operators in Hilbert spaces. We recall that, given
a real Hilbert space $X$, an operator
$A:X\to X$ is said to be a potential operator if there exists a
a G\^ateaux differentiable functional $P$ on $X$ (which is called
a potential of $A$) such that $P'=A$.
\par
\medskip
THEOREM 2.7. - {\it Let $X$ be a real Hilbert
space, and let $A:X\to X$ be a potential operator, with sequentially
weakly upper semicontinuous potential $P$. For each $\rho>0$,
put
$$\varphi(\rho)=\inf_{\|x\|^{2}<\rho}
{{\sup_{\|y\|^{2}\leq \rho}P(y)-P(x)}
\over {\rho-\|x\|^{2}}}\ .$$
Furthermore, set
$$\gamma=\liminf_{\rho\to +\infty}\varphi(\rho)$$
and
$$\delta=\liminf_{\rho\to {0}^{+}}\varphi(\rho)\ .$$
Then, the following conclusions hold:\par
\noindent
$(a)$\hskip 10pt If there is $\rho>0$ such that
$\varphi(\rho)<{{1}\over {2}}$,
then the operator $A$ has a fixed point whose norm is less than
$\sqrt {\rho}$.\par
\noindent
$(b)$\hskip 10pt If $\gamma<{{1}\over {2}}$, then
the following alternative holds: either the functional
$x\to {{1}\over {2}}\|x\|^{2}-P(x)$ has a global minimum, or the set
of all fixed points of $A$ is unbounded.\par
\noindent
$(c)$\hskip 10pt If $\delta<{{1}\over {2}}$, then
 the following alternative holds: either $0$ is 
a local minimum of the functional
$x\to {{1}\over {2}}\|x\|^{2}-P(x)$, or there exists a
sequence of pairwise distinct fixed points of $A$
which weakly converges to $0$.}\par
\smallskip
PROOF. Apply Theorem 2.5 taking $\Phi(x)=-P(x)$ and $\Psi(x)=
\|x\|^{2}$ for all $x\in X$, and observe that $\Psi'(x)=
2x$.\hfill $\bigtriangleup$\par
\medskip
{\it Proof of Theorem A}. We have
$$\liminf_{\rho\to +\infty}\inf_{\|x\|^{2}<\rho}
{{\sup_{\|y\|^{2}\leq \rho}P(y)-P(x)}\over {\rho-\|x\|^{2}}}\leq
\liminf_{\rho\to +\infty}
{{\sup_{\|y\|^{2}\leq \rho}P(y)-P(0)}\over {\rho}}$$
$$=\liminf_{r\to +\infty}{{\sup_{\|x\|\leq r}P(x)}\over
{r^{2}}}<{{1}\over {2}}\ .$$
On the other hand, since $\limsup_{r\to +\infty}
{{\sup_{\|x\|\leq r}P(x)}\over {r^{2}}}>{{1}\over {2}}$, there are a
number $c>{{1}\over {2}}$ and two
sequences $\{r_{n}\}$, $\{x_{n}\}$, with $\lim_{n\to \infty}
r_{n}=+\infty$ and $\|x_{n}\|\leq r_{n}$ for all $n\in {\bf N}$,
such that 
$$P(x_{n})\geq cr_{n}^{2}$$
for all $n\in {\bf N}$. Consequently, we have
$$\lim_{n\to \infty}P(x_{n})=+\infty \eqno {(7)}$$
and
$$P(x_{n})-{{1}\over {2}}\|x_{n}\|^{2}\geq \left (
 c-{{1}\over {2}}\right )
\|x_{n}\|^{2} \eqno {(8)}$$
for all $n\in {\bf N}$. Finally, observe that
the sequence $\{x_{n}\}$ is unbounded.
Indeed, if not, there would be a subsequence $\{x_{n_{k}}\}$ weakly
converging to some $x^{*}$, and so,
 by the sequential weak upper
semicontinuity of $P$, we would have $\limsup_{k\to \infty}
P(x_{n_{k}})\leq P(x^{*})$, against $(7)$. From $(8)$, 
we then infer that the functional $x\to {{1}\over {2}}
\|x\|^{2}-P(x)$ is unbounded below, and so it has no global minima. Now,
the conclusion follows directly from $(b)$ of Theorem 2.7.
\hfill $\bigtriangleup$\par
\medskip
REMARK 2.2. - Another quite recent result on fixed points of potential
operators is provided by Theorem 3.1 of [2]. Among other things,
a major difference
between Theorem A and Theorem 3.1 of [2] is that this latter deals
with weakly continuous potentials.\par
\medskip 
Our final result is a sample of application of Theorem 2.5 to
a class of
nonlinear elliptic equations involving the critical Sobolev exponent.\par
\medskip
THEOREM 2.8. - {\it Let $\Omega\subseteq {\bf R}^{n}$ ($n\geq 3$)
be an open bounded set, with smooth boundary, let $\alpha, \beta\in
L^{{2n}\over {n+2}}(\Omega)$, and let $a, b, c, s, q$ be five real numbers,
with
$b, c>0$, $0<s<1$ and $1<q<{{n+2}\over {n-2}}$. Let $f:{\bf R}\to {\bf R}$
be the function defined by
$$f(\xi)=\cases {\xi^{{n+2}\over {n-2}} & if $\xi\geq 0$\cr &\cr
0 & if $\xi<0 \ .$\cr}$$
Then, there exists $\lambda^{*}>0$ such that, for
each $\lambda\in ]0,\lambda^{*}[$, the problem
$$\cases {-\Delta u=a|u|^{s-1}u+\alpha(x)+\lambda(b|u|^{q-1}u-cf(u)
+\beta(x)) & in $\Omega$\cr & \cr
u_{|\partial \Omega}=0\cr}$$
has at least one weak solution in $W^{1,2}_{0}(\Omega)$.}\par
\smallskip
PROOF. Consider $W^{1,2}_{0}(\Omega)$ endowed with the
norm $\|u\|=(\int_{\Omega}|\nabla u(x)|^{2}dx)^{{1}\over {2}}$. 
Let $g:{\bf R}\to {\bf R}$
be the function defined by
$$g(\xi)=\cases {\xi^{{2n}\over {n-2}} & if $\xi\geq 0$\cr &\cr
0 & if $\xi<0 \ .$\cr}$$
For each $u\in W^{1,2}_{0}(\Omega)$, put
$$\Phi(u)={{c(n-2)}\over {2n}}\int_{\Omega}g(u(x))dx
-{{b}\over {q+1}}\int_{\Omega}|u(x)|^{q+1}dx-\int_{\Omega}
\beta(x)u(x)dx$$
and
$$\Psi(u)={{1}\over {2}}\int_{\Omega}|\nabla u(x)|^{2}dx-
{{a}\over {s+1}}\int_{\Omega}|u(x)|^{s+1}dx-\int_{\Omega}
\alpha(x)u(x)dx\ .$$
By a classical result ([1], Proposition B.10), the
 functionals $\Phi$ and $\Psi$
are G\^ateaux differentiable on $W^{1,2}_{0}(\Omega)$, the weak
solutions of our problem being precisely the critical points
of $\Psi+\lambda\Phi$. Now, observe that, by the Rellich-Kondrachov
theorem, the functionals $u\to \int_{\Omega}|u(x)|^{q+1}dx$
and  $u\to \int_{\Omega}|u(x)|^{s+1}dx$ are sequentially weakly
continuous in $W^{1,2}_{0}(\Omega)$. Moreover, 
the functional $u\to \int_{\Omega}
g(u(x))dx$ is weakly lower semicontinuous  in $W^{1,2}_{0}(\Omega)$,
since it is convex and (strongly) continuous. Hence, the functionals
$\Phi$ and $\Psi$ are sequentially weakly lower semicontinuous
in  $W^{1,2}_{0}(\Omega)$. Finally, observe that, by the Poincar\'e
inequality, since $s<1$, the functional $\Psi$ is coercive. Therefore,
by conclusion $(a)$ of Theorem 2.5, there is a $\mu^{*}>0$ such that,
for each $\mu>\mu^*$, the functional $\Phi+\mu\Psi$ has a critical
point in $W^{1,2}_{0}(\Omega)$. So, we can take $\lambda^{*}=
{{1}\over {\mu^{*}}}$ to get the conclusion.\hfill
$\bigtriangleup$\par
\medskip
REMARK 2.3. - It is worth noticing that, since $b>0$ and
$q>1$,
for each $\lambda>0$,
the functional $\Phi+\lambda\Psi$ in the proof of Theorem 2.8
is unbounded below (and above as well).\par
\bigskip
\bigskip
{\centerline {\bf References}}
\bigskip
\bigskip
\noindent
[1]\hskip 5pt P. H. RABINOWITZ, {\it Minimax methods in critical
point theory with applications to differential equations}, CBMS
Reg. Conf. Ser. in Math., {\bf 65}, Amer. Math. Soc., Providence,
1986.\par
\smallskip
\noindent
[2]\hskip 5pt M. SCHECHTER and K. TINTAREV, {\it Eigenvalues
for semilinear boundary value problems}, Arch. Rational Mech.
Anal., {\bf 113} (1991), 197-208.\par
\bigskip
\bigskip
Department of Mathematics\par
University of Catania\par
Viale A. Doria 6\par
95125 Catania, Italy\par
{\it e-mail address}: ricceri@dipmat.unict.it
\bye